\documentclass[sn-mathphys,Numbered]{sn-jnl}

\usepackage{amssymb, amsmath, amsthm, amscd, mathrsfs}

\usepackage[compatible]{algpseudocode} 
\usepackage{algcompatible}
\usepackage{booktabs}
\usepackage{graphicx}
\usepackage{paralist}
\usepackage{cite}
\usepackage{bigints}
\usepackage{dsfont}
\usepackage{empheq}
\usepackage{amssymb}
\usepackage{cases}
\usepackage{enumitem}
\usepackage{float}
\allowdisplaybreaks
\usepackage{hyperref}

\usepackage{graphicx}%
\usepackage{multirow}%
\usepackage{amsmath,amssymb,amsfonts}%
\usepackage{amsthm}%
\usepackage{mathrsfs}%
\usepackage[title]{appendix}%
\usepackage{xcolor}%
\usepackage{textcomp}%
\usepackage{manyfoot}%
\usepackage{booktabs}%
\usepackage{algorithm}%
\usepackage{algorithmicx}%
\usepackage{algpseudocode}%
\usepackage{listings}%

\theoremstyle{plain}
\newtheorem{teo}{Theorem}[section]

\newtheorem{lemma}[teo]{Lemma}

\theoremstyle{definition}
\newtheorem{defi}[teo]{Definition}

\begin{document}

\title[Numerical impulse controllability for parabolic equations]
      {Numerical impulse controllability for parabolic equations by a penalized HUM approach
}


\author[1]{\fnm{Salah-Eddine} \sur{Chorfi}}\email{s.chorfi@uca.ac.ma}

\author[2]{\fnm{Ghita} \sur{El Guermai}}\email{ghita.el.guermai@gmail.com}

\author[3]{\fnm{Lahcen} \sur{Maniar}}\email{maniar@uca.ac.ma}

\author[4]{\fnm{Walid} \sur{Zouhair}}\email{walid.zouhair.fssm@gmail.com}

\affil[1,2,3]{\orgdiv{LMDP, UMMISCO (IRD-UPMC)}, \orgname{Cadi Ayyad University, Faculty of Sciences Semlalia}, \orgaddress{\city{Marrakesh}, \postcode{B.P. 2390}, \country{Morocco}}}

\affil[4]{\orgdiv{Department of mathematics}, \orgname{Ibn Zohr University, Faculty of Applied Sciences Ait Melloul}, \orgaddress{\street{Route Nationale N10}, \city{Azrou}, \postcode{B.P. 6146}, \country{Morocco}}}



\abstract{This work presents a comparative study to numerically compute impulse approximate controls for parabolic equations with various boundary conditions. Theoretical controllability results have been recently investigated using a logarithmic convexity estimate at a single time based on a Carleman commutator approach. We propose a numerical algorithm for computing the impulse controls with minimal $L^2$-norms by adapting a penalized Hilbert Uniqueness Method (HUM) combined with a Conjugate Gradient (CG) method. We consider static boundary conditions (Dirichlet and Neumann) and dynamic boundary conditions. Some numerical experiments based on our developed algorithm are given to validate and compare the theoretical impulse controllability results.}


\maketitle

\section{Introduction and main results}
Parabolic equations, where the heat equation is the prototype, constitute a class of Partial Differential Equations (PDEs) that describe the evolution of physical quantities over time and space. The heat equation is particularly important in the study of heat transfer and diffusion processes.

Impulsive systems in the context of PDEs refer to systems whose behavior changes abruptly or impulsively at certain points in space or time. These impulsive changes can be modeled mathematically using PDEs with discontinuities or Dirac delta functions, which are used to present concentrated impulses at specific spatial or temporal locations. Impulsive systems are encountered in various fields, from biological models to fluid dynamics as well as economics, among others. They manifest as sudden boundary conditions in switched control inputs or shock waves in compressible flows. Controlling impulsive systems via impulse controls can be challenging (due to the presence of delta functions) requiring specific techniques such as logarithmic convexity, Carleman commutator approach, and optimal impulse control theory, see e.g., \citep{XLSS, BCKDP, ABWZ}.

Impulsive controllability is a concept in control theory that deals with the ability to control a dynamical system by applying control inputs at specific discrete instants or time intervals, often referred to as ``impulse times" (or intervals). In impulsive controllability, the control actions are not continuously applied but occur at distinct time points. Impulsive systems can model situations where control inputs change abruptly. They make it possible to handle systems in situations where continuous control may not be feasible or practical. For instance, in switched control systems, the control input can change instantaneously at a specific instant. In this context, the impulse approximate controllability was studied for a linear heat equation with homogeneous Dirichlet and Neumann boundary conditions in \citep{pkm, RBKDP}, using a new strategy combining the logarithmic convexity method and the Carleman commutator approach. In \citep{Buffephung}, the authors have established a Lebeau-Robbiano-type spectral inequality for a degenerate one-dimensional elliptic operator with application to impulse control and finite-time stabilization. It should be pointed out that this method is a new approach to steer the solution to zero using impulse control as a stabilizer in finite time. Recently, in \citep{CGMZ'21, CGMZ'23, CGMZ'22} the authors have established new results of impulse controllability for a general type of dynamic boundary conditions which introduce mathematical issues that require sophisticated estimates due to the boundary terms. We refer to the seminal paper \citep{MMS'17} for more details on the non-impulsive control case. 

Boundary conditions play a crucial role in solving PDEs as they have a significant effect on the behavior of the solutions. They describe the interaction of the system with its surroundings. For example, in heat transfer problems, the temperature on the boundary may represent an insulated or constant temperature boundary, reflecting the physical properties of the system. The choice of appropriate boundary conditions is a fundamental step in the analysis of PDEs in various fields of science and engineering. From a numerical perspective, boundary conditions are fundamental for ensuring the accuracy, stability, and convergence of numerical solutions to PDEs. They guide how the spatial domain should be discretized. In this study, we propose an algorithm designed for the numerical computation of impulse controls with minimal energy. This approach involves an adaptation of the penalized HUM and a CG method to the impulsive case. For further information, we recommend the book \citep{GL'08} and the paper \citep{Bo'13}. Our investigation encompasses both static boundary conditions (Dirichlet and Neumann) as well as dynamic boundary conditions. To validate and compare theoretical findings regarding impulse controllability, we conduct several numerical experiments using the algorithm we have developed. Finally, it should be emphasized that the numerical computation of impulse controls has not been considered before for static boundary conditions. We refer to the section ``Future works" of the thesis \citep{Vo'18}.

The remainder of this paper is structured as follows: in Section \ref{sec2}, we provide a review of various results related to impulse controllability for the heat equation with different boundary conditions. Section \ref{sec3} is devoted to the algorithm for computing impulse optimal controls, accompanied by numerical simulations for illustration. Finally, we conclude with a comparative analysis of numerical outcomes across various boundary conditions.

\section{Preliminary results}\label{sec2}
Let $\Omega\subset \mathbb{R}^n$ be a bounded domain with smooth boundary $\Gamma:=\partial\Omega$. Let $T>0$ be an arbitrary control time and $\tau \in (0, T)$ be an arbitrary fixed impulsion time. We consider the following impulse-controlled system
\begin{empheq}[left = \empheqlbrace]{alignat=2} \label{eq1.1}
\begin{aligned}
&\partial_{t} \Psi- \mathbf{A} \Psi=0, && \qquad\text { in } (0, T) \backslash\{\tau\},\\
&\psi(\cdot, \tau)=\psi\left(\cdot, \tau^{-}\right)+\mathds{1}_{\omega} h(\cdot,\tau), && \qquad\text { in } \Omega,\\
&\Psi(0) = \Psi^0 , \\
\end{aligned}
\end{empheq}
where $\psi(\cdot,\tau^{-})$ denotes the left limit of the function $\psi$ at time $\tau$, the control region $\omega\Subset \Omega$ is a nonempty open subset, $\mathds{1}_\omega$ stands for the characteristic function of $\omega$,  $h(\cdot,\tau)$ is an impulsive control acting at the impulse instant $\tau$. The notation $\mathbf{A}$ designates a linear operator on an $L^2$-space with norm $\|\cdot\|$, and the state $\Psi$ might be a couple $(\psi, \psi_\Gamma)$, depending on the type of boundary conditions (see the next subsections).

If the operator $\mathbf{A}$ generates a $C_0$-semigroup. Then, for every initial datum $\Psi_{0}$, the system \eqref{eq1.1} has a unique mild solution given by
$$
\Psi(t) = \mathrm{e}^{t\mathbf{A}} \Psi_{0} + \mathds{1}_{\{t\geq \tau \}}(t)\, \mathrm{e}^{(t-\tau)\mathbf{A}} (\mathds{1}_{\omega} h(\tau),0), \qquad t\in (0,T).
$$
\begin{defi}
System \eqref{eq1.1} is null approximate impulse controllable at time $T$ if for any $\varepsilon > 0$ and
any $\Psi^0$, there exists a control function $h \in L^2(\omega),$  such that the associated state at final time satisfies
\begin{equation*}
\|\Psi(\cdot, T)\| \leq \varepsilon\left\Vert \Psi^0\right\Vert .
\end{equation*}
\end{defi}
This means that for every $\varepsilon >0$ and every initial datum $\Psi^0$, the set
\begin{equation*}
\mathcal{R}_{T, \Psi^{0}, \varepsilon} :=\left\{h \in L^{2}(\omega): \text { the solution of }\eqref{eq1.1}\text { satisfies }\left\Vert\Psi(\cdot, T)\right\Vert \leq \varepsilon\left\Vert \Psi^{0}\right\Vert\right\},
\end{equation*}
is nonempty; which leads to the definition of the cost of null approximate impulse controllability.
\begin{defi}
The quantity  $$K(T,\varepsilon):=\sup_{\left\|\Psi^0\right\|=1} \inf_{h \in \mathcal{R}_{T, \Psi^{0}, \varepsilon}} \|h\|_{L^{2}(\omega)}
$$ is called the cost of null
approximate impulse controllability at time $T.$
\end{defi}

\subsection{Dirichlet case}
In this subsection, we recall the impulsive controllability result for the heat equation with the Dirichlet boundary condition:
\begin{empheq}[left = \empheqlbrace]{alignat=2} \label{d1.1}
\begin{aligned}
&\partial_{t} \psi-\Delta \psi=0, && \qquad\text { in } \Omega \times(0, T) \backslash\{\tau\},\\
&\psi(\cdot, \tau)=\psi\left(\cdot, \tau^{-}\right)+\mathds{1}_{\omega} h(\cdot,\tau), && \qquad\text { in } \Omega,\\
& \psi = 0, &&\qquad\text{ on } \Gamma \times(0, T) , \\
&\psi(\cdot, 0) = \psi^{0}, && \qquad \text{ on } \Omega.
\end{aligned}
\end{empheq}
To obtain the null approximate impulse controllability of the above equation, the key ingredient is the following logarithmic convexity estimate:
\begin{lemma}[\citep{pkm}]
For any $T > 0$ and any $\omega$ nonempty open subset of $\Omega$,
\begin{equation}\label{e1.2}
\|u(\cdot, T)\|_{L^2(\Omega)} \leq \mathrm{e}^{C \frac{K}{T}}\|u(\cdot, T)\|_{L^{2}(\omega)}^{\beta}\|u(\cdot, 0)\|^{1-\beta}_{L^2(\Omega)},
\end{equation}
where $\beta \in (0,1)$, $C, K >0$ are constants, and $u$ is the solution of the following homogeneous system
\begin{empheq}[left = \empheqlbrace]{alignat=2}
\begin{aligned}
&\partial_{t} u-\Delta u=0, && \qquad\text { in } \Omega \times(0, T), \\
& u  = 0, &&\qquad\text{ on } \Gamma \times(0, T) , \\
& u(\cdot, 0)=u^{0}, && \qquad \text{ in } \Omega.
\end{aligned}
\end{empheq}
\end{lemma}
The previous lemma reflects an observability estimate at a single instant of time. This estimate has been proven using the weight function
\[
\Phi (x,t)=\frac{-\left\lvert x-x_{0}\right\lvert^2}{4(T-t+\rho)}, \qquad (x,t) \in \overline{\Omega}\times (0,T),
\]  
where $x_0\in \omega$ and $\rho>0$ is suitably chosen.

Consequently, the following result on impulse controllability of the equation \eqref{d1.1} was established:
\begin{teo}[\citep{pkdgwyx, Vo'17}]
The heat equation \eqref{d1.1} is null approximate impulse controllable at time $T$. Moreover, we have the following upper bound for the control cost 
$$K_1(T, \varepsilon) \leq \frac{M_{1} \mathrm{e}^{\frac{M_{2}}{T-\tau}}}{\varepsilon^{\delta}},$$
where $M_1$, $M_2$ and $\delta$ are positive constants depending on $\Omega$ and $\omega$.
\end{teo}

\subsection{Neumann case}
Here, we recall the impulsive controllability result for the heat equation with the Neumann boundary condition:
\begin{empheq}[left = \empheqlbrace]{alignat=2} \label{d1.5}
\begin{aligned}
&\partial_{t} \psi-\Delta \psi=0, && \qquad\text { in } \Omega \times(0, T) \backslash\{\tau\},\\
&\psi(\cdot, \tau)=\psi\left(\cdot, \tau^{-}\right)+\mathds{1}_{\omega} h(\cdot,\tau), && \qquad\text { in } \Omega,\\
& \partial_{\nu}\psi = 0, &&\qquad\text{ on } \Gamma \times(0, T) , \\
&\psi(\cdot, 0) = \psi^{0}, && \qquad \text{ in } \Omega,
\end{aligned}
\end{empheq}
where $\nu$ is the unit outward normal vector to $\Gamma$, and $\partial_\nu \psi$ denotes the normal derivative. As before, the key ingredient is the following logarithmic convexity estimate:
\begin{lemma}[\citep{RBKDP}]
For any $T > 0$ and any $\omega$ nonempty open subset of $\Omega$,
\begin{equation*}
\|u(\cdot, T)\|_{L^2(\Omega)} \leq\left(e^{C\left(1+\frac{1}{T}\right)}\|u(\cdot, T)\|_{L^2(\omega)}\right)^\beta\|u(\cdot, 0)\|_{L^2(\Omega)}^{1-\beta} .
\end{equation*}
where $\beta \in (0,1)$, $C >0$ are constants only depending on $\Omega$ and $\omega$, and $u$ is the solution of the homogeneous  system
\begin{empheq}[left = \empheqlbrace]{alignat=2}
\begin{aligned}
&\partial_{t} u-\Delta u=0, && \qquad\text { in } \Omega \times(0, T), \\
& \partial_{\nu}u = 0, &&\qquad\text{ on } \Gamma \times(0, T) , \\
& u(\cdot, 0)=u^{0}, && \qquad \text{ in } \Omega.
\end{aligned}
\end{empheq}
\end{lemma}
This result has been recently extended to a general parabolic equation with variable diffusion and drift coefficients in \citep{Du22}.

Being different from the Dirichlet case, the above lemma has been established by introducing a small parameter $s\in (0,1)$ in the weight function
\[
\Phi_s(x,t)=\frac{-s\left\lvert x-x_{0}\right\lvert^2}{4(T-t+\rho)}, \qquad (x,t) \in \overline{\Omega}\times (0,T),
\]  
where $x_0\in \omega$ and $\rho>0$ is suitably chosen.

Then, one can prove the following impulse controllability for the equation \eqref{d1.5}:
\begin{teo}[\citep{RBKDP}]
The heat equation \eqref{d1.5} is null approximate impulse controllable at time $T$. Moreover, we have the following upper bound for the control cost 
$$K_2(T, \varepsilon) \leq \frac{N_{1} \mathrm{e}^{\frac{N_{2}}{T-\tau}}}{\varepsilon^{\sigma}},$$
where $N_1$, $N_2$ and $\sigma$ are positive constants depending on $\Omega$ and $\omega$.
\end{teo}

\subsection{Dynamic case}
Now, we consider the following heat equation with dynamic boundary conditions
\begin{empheq}[left = \empheqlbrace]{alignat=2} \label{1.1}
\begin{aligned}
&\partial_{t} \psi-\Delta \psi=0, && \qquad\text { in } \Omega \times(0, T) \backslash\{\tau\},\\
&\psi(\cdot, \tau)=\psi\left(\cdot, \tau^{-}\right)+\mathds{1}_{\omega} h(\cdot,\tau), && \qquad\text { in } \Omega,\\
&\partial_{t}\psi_{\Gamma} - \Delta_{\Gamma} \psi_{\Gamma} + \partial_{\nu}\psi =0, && \qquad\text { on } \Gamma \times(0, T)\backslash\{\tau\}, \\
&\psi_{\Gamma}(\cdot, \tau)=\psi_{\Gamma}\left(\cdot, \tau^{-}\right), && \qquad\text { on } \Gamma,\\
& \psi_{\Gamma}(x,t) = \psi_{|\Gamma}(x,t), &&\qquad\text{ on } \Gamma \times(0, T) , \\
& \left(\psi(\cdot, 0),\psi_{\Gamma}(\cdot, 0)\right)=\left(\psi^{0},\psi^{0}_{\Gamma}\right), && \qquad \text{ on } \Omega\times\Gamma,
\end{aligned}
\end{empheq}
where $\left(\psi^{0},\psi^{0}_{\Gamma}\right)\in L^2(\Omega)\times L^2(\Gamma)$ denotes the initial condition. Again, the key result is the logarithmic convexity estimate.
\begin{lemma}[\citep{CGMZ'21, CGMZ'23}]
For any $T > 0$ and any $\omega$ nonempty open subset of $\Omega$, the following estimate holds
\begin{equation}\label{1.3}
\|U(\cdot, T)\|_{L^2(\Omega)\times L^2(\Gamma)} \leq\left(\mu \mathrm{e}^{\frac{K}{T}}\|u(\cdot, T)\|_{L^{2}(\omega)}\right)^{\beta}\|U(\cdot, 0)\|^{1-\beta}_{L^2(\Omega)\times L^2(\Gamma)},
\end{equation}
where $\mu, K >0$, $\beta \in (0,1)$ are constants, and $U=\left(u,u_{\Gamma}\right)$ is the solution of the following homogeneous system
\begin{empheq}[left = \empheqlbrace]{alignat=2}\label{e1.3}
\begin{aligned}
&\partial_{t} u-\Delta u=0, && \qquad\text { in } \Omega \times(0, T), \\
&\partial_{t}u_{\Gamma} - \Delta_{\Gamma} u_{\Gamma} + \partial_{\nu}u =0, && \qquad\text { on } \Gamma \times(0, T), \\
& u_{\Gamma}(x,t) = u_{|\Gamma}(x,t), &&\qquad\text{ on } \Gamma \times(0, T) , \\
& \left(u(\cdot, 0),u_{\Gamma}(\cdot, 0)\right)=\left(u^{0},u^{0}_{\Gamma}\right), && \qquad \text{ on } \Omega\times\Gamma.
\end{aligned}
\end{empheq}
\end{lemma}
In this dynamic case, several new boundary terms occur and should be absorbed. This has been done thanks to the small parameter $s$ introduced in the weight function $\Phi_s$ inspired by the Neumman case.

Consequently, we obtained the following impulse controllability result:
\begin{teo}[\citep{CGMZ'21}]\label{thm1.1}
The system \eqref{1.1} is null approximate impulse controllable at any time $T > 0$. Moreover, we have the following upper bound for the control cost 
$$K_3(T, \varepsilon) \leq \frac{L_{1} \mathrm{e}^{\frac{L_{2}}{T-\tau}}}{\varepsilon^{\kappa}},$$
where $L_1$, $L_2$ and $\kappa$ are positive constants depending on $\Omega$ and $\omega$.
\end{teo}

\section{Algorithm for calculating HUM impulse controls}\label{sec3}
In this section, we propose a numerical algorithm designed for calculating the HUM impulse controls. This method employs a penalized HUM approach along with a CG algorithm. We refer to \citep{GL'08} and \citep{Bo'13} for more details on such a method.

\subsection*{Notations}
We introduce the following notations to encapsulate various boundary conditions and give a general algorithm:
$$\mathbb{L}^2:=\begin{cases}
    L^2(\Omega), \hspace{2cm} (\text{Dirichlet and Neumann cases}),\\
    L^2(\Omega)\times L^2(\Gamma), \qquad (\text{Dynamic case}),
\end{cases}$$
with the inner product
$$\langle \cdot, \cdot\rangle:=\begin{cases}
    \langle \cdot, \cdot\rangle_{L^2(\Omega)}, \hspace{2.5cm} (\text{Dirichlet and Neumann cases}),\\
    \langle \cdot, \cdot\rangle_{L^2(\Omega)}+\langle \cdot, \cdot\rangle_{L^2(\Gamma)}, \qquad (\text{Dynamic case})
\end{cases}$$
and the norm
$$\|\cdot\|:=\begin{cases}
    \|\cdot\|_{L^2(\Omega)}, \hspace{1.7cm} (\text{Dirichlet and Neumann cases}),\\
    \|\cdot\|_{L^2(\Omega)\times L^2(\Gamma)}, \qquad (\text{Dynamic case}).
\end{cases}$$
Any capital letter as $\vartheta$ will stand for the couple $(\upsilon, \upsilon_\Gamma) \in L^2(\Omega)\times L^2(\Gamma)$. In particular, we will identify $(\upsilon, \upsilon_\Gamma) \in L^2(\Omega)\times L^2(\Gamma)$ with $\upsilon \in L^2(\Omega)$ in Dirichlet and Neumann cases. We denote by $\mathbf{BC}$ one of the boundary conditions: Dirichlet condition, Neumann condition, or Dynamic condition. In each case, the operator $\mathbf{A}$ stands for the governing linear operator, and $\mathrm{e}^{t \mathbf{A}}$ designates its associated $C_0$-semigroup on $\mathbb{L}^2$.

\subsection{The HUM impulse controls}
Let $\varepsilon>0$ be fixed and let $\Psi^0$ be an initial datum to be controlled. Without loss of generality, we may assume that $\|\Psi^0\|=1$. We define the cost functional $J_{\varepsilon}: \mathbb{L}^2 \rightarrow \mathbb{R}$ by
\begin{equation*}
J_{\varepsilon}\left(\vartheta^{0}\right)=\frac{1}{2}\|\upsilon(\cdot, T-\tau)\|_{L^{2}(\omega)}^{2}+\frac{\varepsilon}{2}\left\|\vartheta^{0}\right\|^2 + \left\langle \Psi^{0}, \vartheta(\cdot,T) \right\rangle,
\end{equation*}
where $\vartheta$ is the solution of the homogeneous heat equation with $\mathbf{BC}$ corresponding to $\vartheta^{0}$. Note that the functional $J_{\varepsilon}$ is strictly convex, of class $C^1$, and coercive, i.e., $J_{\varepsilon}\left(\vartheta^{0}\right) \to \infty$ as $\|\vartheta^{0}\|\to \infty$. Then the unique minimizer $\tilde{\vartheta}^{0}_\varepsilon \in \mathbb{L}^2$ of $J_{\varepsilon}$ is characterized by the Euler-Lagrange equation
\begin{equation}\label{Eq2}
\int_{\omega} \tilde{\upsilon}_\varepsilon(x, T-\tau) z(x, T-\tau) \mathrm{d} x+\varepsilon \left\langle \tilde{\vartheta}^{0}_\varepsilon,Z^{0} \right\rangle+\left\langle \Psi^{0}, Z(\cdot,T) \right\rangle = 0
\end{equation}
for all $Z^0 \in \mathbb{L}^2$, where $\tilde{\vartheta}_\varepsilon$ and $Z$ are respectively the solutions of the homogeneous heat equation with $\mathbf{BC}$ corresponding to $\tilde{\vartheta}^{0}_\varepsilon$ and $Z^{0}$. We introduce the control operator $\mathcal{B}\colon \mathbb{L}^2 \rightarrow \mathbb{L}^2$ defined by
$$\mathcal{B} \vartheta=(\mathds{1}_\omega \upsilon,0),$$
and we consider the non-negative symmetric operator (the Gramian operator)
$$\Lambda_\tau \colon \mathbb{L}^2 \rightarrow \mathbb{L}^2,$$
given by
$$\Lambda_\tau \varrho=\mathrm{e}^{(T-\tau)\mathbf{A}} \mathcal{B} \,\mathrm{e}^{(T-\tau)\mathbf{A}} \varrho.$$
Thus, the HUM impulse control is given by
$$\widehat{h}=\mathcal{B} \,\mathrm{e}^{(T-\tau)\mathbf{A}} \tilde{\vartheta}^{0}_\varepsilon,$$
and the identity \eqref{Eq2} can be rewritten as
\begin{equation}
    \left(\Lambda_\tau +\varepsilon \mathbf{I}_{\mathbb{L}^2}\right) \tilde{\vartheta}^{0}_\varepsilon = -\mathrm{e}^{T\mathbf{A}} \Psi^0,
\end{equation}
where $\mathbf{I}_{\mathbb{L}^2}$ denotes the identity operator. To resolve the above operator equation, we propose the following CG algorithm.
\bigskip

\begin{small}
\begin{algorithm}[H]\label{alg1}
 Set $k=0$ and choose an initial guess $\mathbf{f}_0=\left(f_0,f_{0,\Gamma}\right) \in \mathbb{L}^2$.\;
 Solve the problem
 \begin{empheq}[left = \empheqlbrace]{alignat=2}
\begin{aligned}
&\partial_{t} p_0(x,t)-\Delta p_0(x,t)=0, && \quad (x,t)\in \Omega_T, \\
&\mathbf{BC}, && \quad t\in (0, T), \\
&\left(p_0(x,0),p_{0,\Gamma}(x, 0)\right)=\mathbf{f}_0(x), && \quad x\in \Omega, \notag
\end{aligned}
\end{empheq}
and set $\mathbf{u}_0(x)=\mathcal{B}\, \mathbf{p}_0(T-\tau,x)$. \; 
 Solve the problem
 \begin{empheq}[left = \empheqlbrace]{alignat=2}
\begin{aligned}
&\partial_{t} y_0(x,t)-\Delta y_0(x,t)=0, && \, (x,t)\in \Omega_T , \\
&\mathbf{BC}, && \, t\in (0, T), \\
&\left(y_0(x,0),y_{0,\Gamma}(x, 0)\right)=\mathbf{u}_0(x), && \, x\in \Omega, \notag
\end{aligned}
\end{empheq}
compute $\mathbf{g}_0=\varepsilon \mathbf{f}_0 + Y_0(T-\tau) + \mathrm{e}^{T\mathbf{A}} \Psi^0$ and set $\mathbf{w}_0=\mathbf{g}_0$.
\;
 For $k=1,2,\ldots,$ until convergence,
 solve the problem
 \begin{empheq}[left = \empheqlbrace]{alignat=2}
\begin{aligned}
&\partial_{t} p_k(x,t)-\Delta p_k(x,t)=0, && \, (x,t)\in \Omega_T , \\
&\mathbf{BC}, && \, t\in (0, T), \\
&\left(p_k(x,0),p_{k,\Gamma}(x, 0)\right)=\mathbf{w}_{k-1}(x), && \, x\in \Omega, \notag
\end{aligned}
\end{empheq}
and set $\mathbf{u}_k(x)=\mathcal{B}\, \mathbf{p}_k(T-\tau,x)$.\;
Solve the problem
 \begin{empheq}[left = \empheqlbrace]{alignat=2}\label{1.3}
\begin{aligned}
&\partial_{t} y_k(x,t)-\Delta y_k(x,t)=0, && \, (x,t)\in \Omega_T , \\
&\mathbf{BC}, && \, t\in (0, T), \\ 
&\left(y_k(x,0),y_{k,\Gamma}(x, 0)\right)=\mathbf{u}_k(x), && \, x\in \Omega, \notag
\end{aligned}
\end{empheq}
and compute $$\bar{\mathbf{g}}_k=\varepsilon \mathbf{w}_{k-1} + Y_k(T-\tau) \qquad \text{ and } \qquad \rho_k=\dfrac{\|\mathbf{g}_{k-1}\|^2}{\langle \bar{\mathbf{g}}_k, \mathbf{w}_{k-1} \rangle},$$
then
$$\mathbf{f}_k=\mathbf{f}_{k-1}-\rho_k \mathbf{w}_{k-1} \qquad \text{ and } \qquad \mathbf{g}_k=\mathbf{g}_{k-1}-\rho_k \bar{\mathbf{g}}_k.$$\;
 \textbf{If} $\dfrac{\|\mathbf{g}_{k}\|}{\|\mathbf{g}_{0}\|} \le tol$, stop the algorithm, set $\mathbf{g}=\mathbf{f}_{k}$ and solve the problem
 \begin{empheq}[left = \empheqlbrace]{alignat=2}\label{1.3}
\begin{aligned}
&\partial_{t} p_k(x,t)-\Delta p_k(x,t)=0, && \, (x,t)\in \Omega_T , \\
&\mathbf{BC}, && \, t\in (0, T), \\
&\left(p_k(x,0),p_{k,\Gamma}(x, 0)\right)=\mathbf{g}(x), && \, x\in \Omega, \notag
\end{aligned}
\end{empheq}
and set $\mathbf{u}_k(x)=\mathcal{B}\, \mathbf{p}_k(T-\tau,x)$. \newline
\textbf{Else} compute $$\gamma_k=\dfrac{\|\mathbf{g}_{k}\|^2}{\|\mathbf{g}_{k-1}\|^2} \qquad \text{ and then } \qquad \mathbf{w}_k=\mathbf{g}_{k}+\gamma_k \mathbf{w}_{k-1}.$$
 \caption{HUM with CG Algorithm}
\end{algorithm}
\end{small}

\subsection{Numerical experiments}
Now, we conduct several numerical tests to demonstrate the theoretical findings and to highlight the effectiveness of the above CG algorithm.

In all main numerical experiments, we will take the following values
$$T=0.02, \quad \tau=0.01, \quad \Omega=(0,1), \quad \omega=(0.3, 0.7) \Subset (0,1),$$
and the initial datum to be controlled is given by
$$\psi_0(x)=\sqrt{2} \sin(\pi x), \qquad x \in [0,1].$$

We employ the method of lines to numerically solve diverse parabolic equations subject to different boundary conditions in Algorithm \ref{alg1}. In this approach, we use the uniform spatial grid given by $x_j=j \Delta x$ for  $j=\overline{0, N_x}$, with $\Delta x=\dfrac{1}{N_x}$. Next, we denote by $u_j(t):=u(t,x_j)$. The second-order derivative of $u$ is approximated by
$$u_{xx}(t,x_j) \approx \frac{u_{j-1}(t)-2 u_j(t) + u_{j+1}(t)}{(\Delta x)^2}, \quad j=\overline{1, N_x-1}.$$
The first-order derivatives on the boundary are approximated by
\begin{align*}
    u_x(t,0) &\approx \frac{u_1(t)-u_0(t)}{\Delta x}\\
    u_x(t,1) &\approx \frac{u_{N_x}(t)-u_{N_x-1}(t)}{\Delta x}.
\end{align*}
Thus, it suffices to resolve the resulting system of ordinary differential equations.

For our computations, we take $N_x=25$ for the spatial mesh parameter. The initial guess in the algorithm is taken as $\mathbf{f}_0=0$. We also choose $\varepsilon=10^{-2}$ and the stopping parameter $tol=10^{-3}$ for the plots.

\subsection{Dirichlet case}
We plot the uncontrolled and the controlled solutions.
\begin{figure}[H]
\centering
\begin{minipage}{0.45\textwidth}
\includegraphics[scale=0.5]{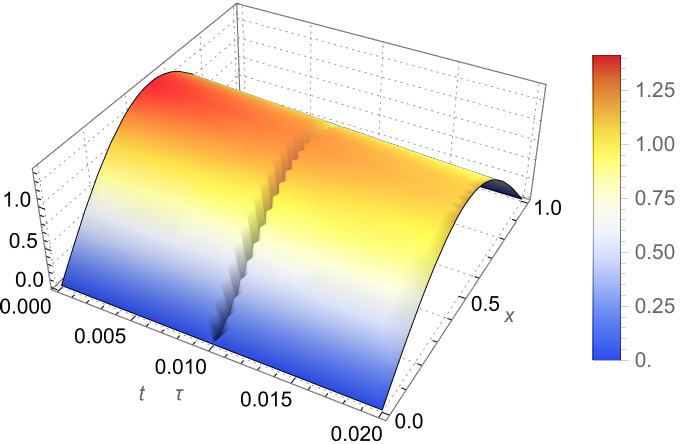}
\caption{The uncontrolled solution in Dirichlet case.}
\end{minipage}\hfill
\begin{minipage}{0.45\textwidth}
\centering
\includegraphics[scale=0.5]{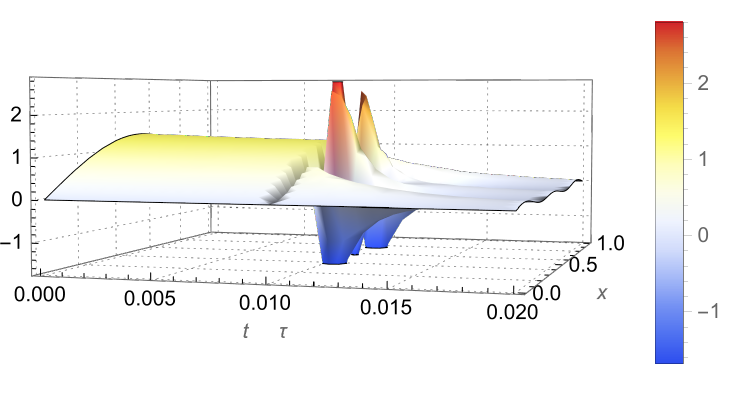}
\caption{The controlled solution in Dirichlet case.}
\label{fig1dc}
\end{minipage}
\end{figure}
The algorithm stops at the iteration number $k_*=10$.

\begin{figure}[H]
\centering
\begin{minipage}{0.45\textwidth}
\includegraphics[scale=0.4]{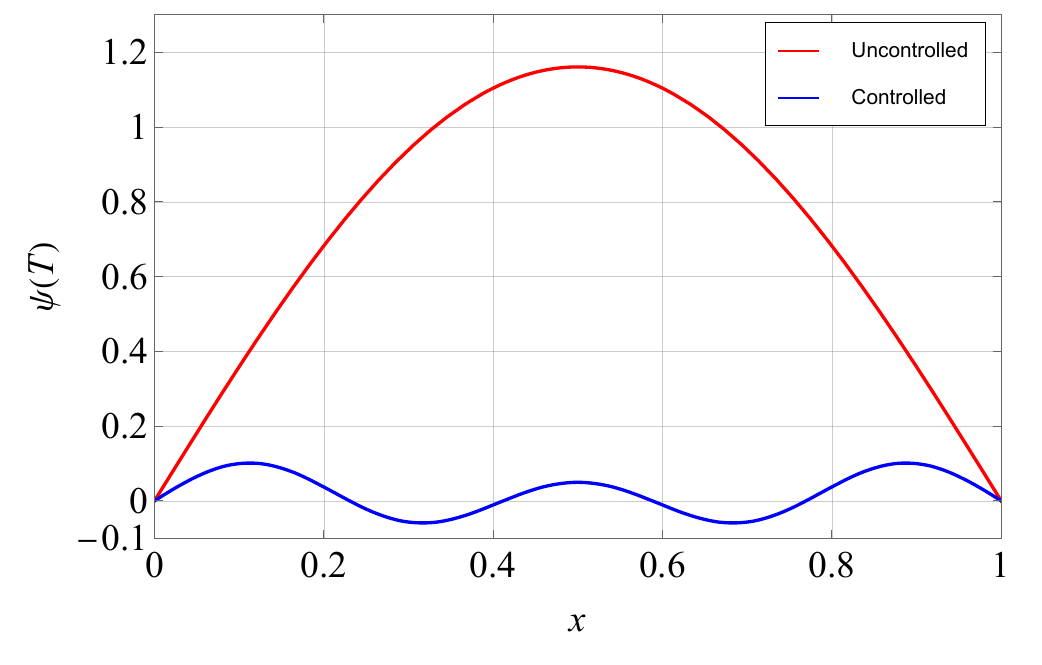}
\caption{The final state for uncontrolled and controlled solutions in Dirichlet case.}
\end{minipage}\hfill
\begin{minipage}{0.45\textwidth}
\centering
\includegraphics[scale=0.4]{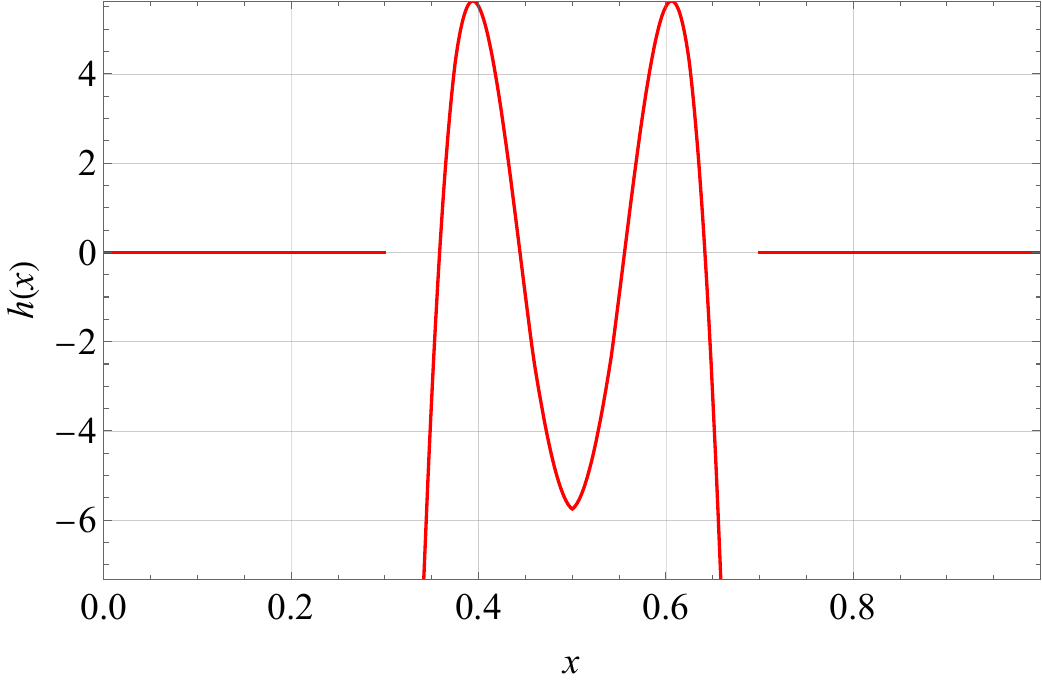}
\caption{The computed control $h$ in Dirichlet case.}
\end{minipage}
\end{figure}
\bigskip

\begin{table}[ht]
\caption{Numerical results for $T=0.02$, $\tau=0.01$ and $tol=10^{-3}$ for Dirichlet condition.} \label{t1}
\centering\setlength\arraycolsep{13pt}
\begin{tabular}{cccc}
\hline\\[-3mm]
 $\varepsilon$  & $10^{-1}$ & $10^{-2}$ & $10^{-3}$ \\
\hline\\[-3mm]
 $N_{\text{iter}}$ & $4$ & $10$ & $20$ \\
\hline\\[-3mm]
 $\|\psi_{\mathrm{D}}(T)\|$ & $1.148\times 10^{-1}$ & $5.63\times 10^{-2}$ & $1.73\times 10^{-2}$ \\
\hline\\[-3mm]
 $\|h_{\mathrm{D}}\|_{L^2(\omega)}$ & $1.684$ & $7.3014$ & $28.6994$ \\
\hline\\[-3mm]
\end{tabular}
\end{table}

\subsection{Neumann case}
Next, we plot the uncontrolled and the controlled solutions.
\begin{figure}[H]
\centering
\begin{minipage}{0.45\textwidth}
\includegraphics[scale=0.5]{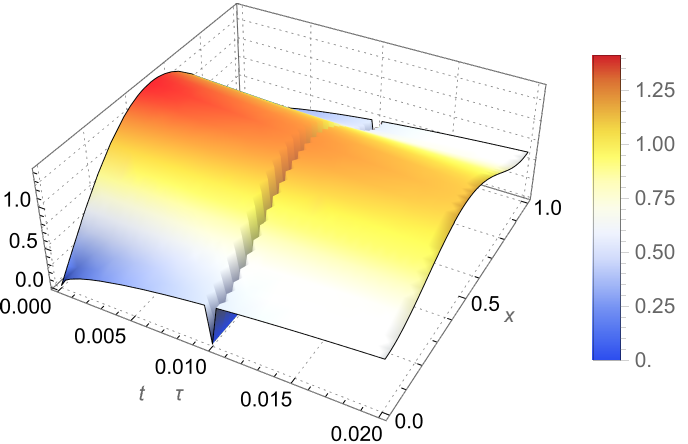}
\caption{The uncontrolled solution in Neumann case.}
\end{minipage}\hfill
\begin{minipage}{0.45\textwidth}
\centering
\includegraphics[scale=0.5]{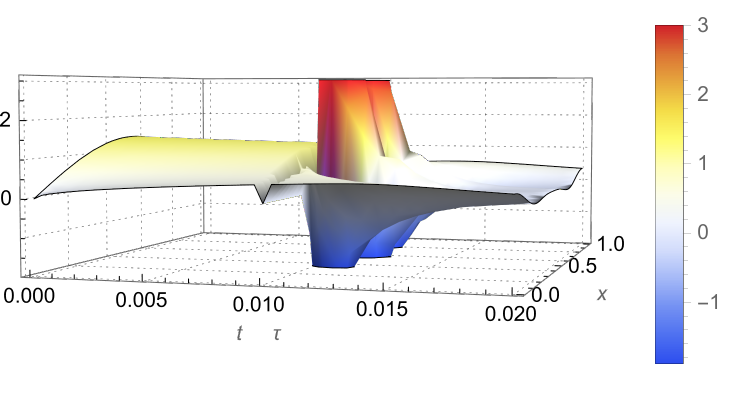}
\caption{The controlled solution in Neumann case.}
\label{fig1nc}
\end{minipage}
\end{figure}
The algorithm stops at the iteration number $k_*=29$. 

\begin{figure}[H]
\centering
\begin{minipage}{0.45\textwidth}
\includegraphics[scale=0.4]{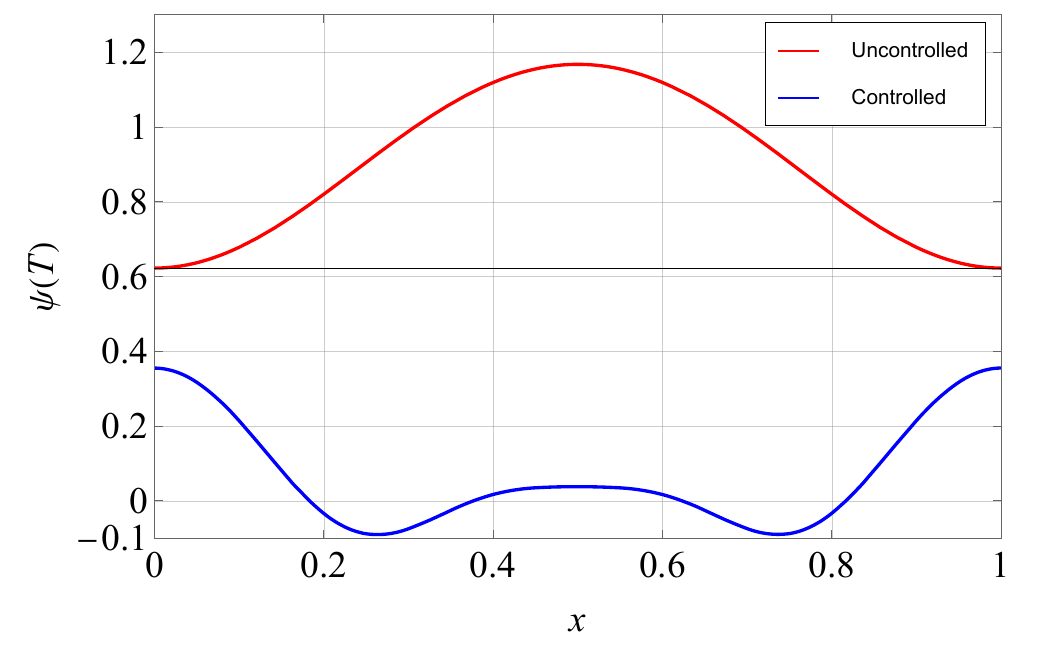}
\caption{The final state for uncontrolled and controlled solutions in Neumann case.}
\end{minipage}\hfill
\begin{minipage}{0.45\textwidth}
\centering
\includegraphics[scale=0.4]{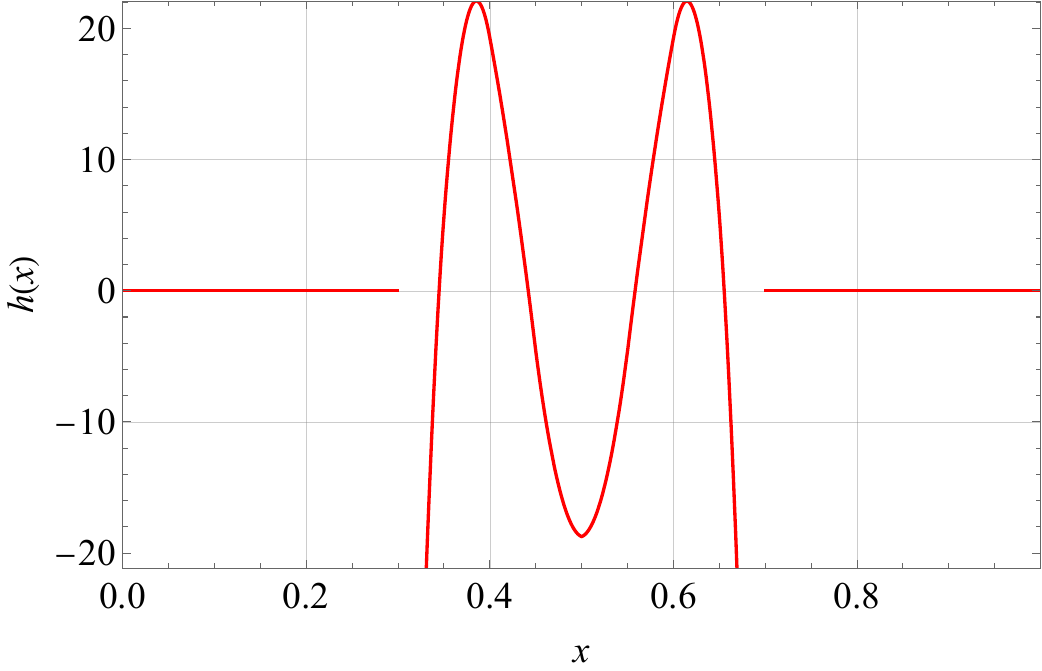}
\caption{The computed control $h$ in Neumann case.}
\end{minipage}
\end{figure}
\bigskip

\begin{table}[ht]
\caption{Numerical results for $T=0.02$, $\tau=0.01$ and $tol=10^{-3}$ for Neumann condition.}
\label{t2}
\centering\setlength\arraycolsep{13pt}
\begin{tabular}{cccc}
\hline\\[-3mm]
 $\varepsilon$  & $10^{-1}$ & $10^{-2}$ & $10^{-3}$ \\
\hline\\[-3mm]
 $N_{\text{iter}}$ & $4$ & $29$ & $100$ \\
\hline\\[-3mm]
 $\|\psi_{\mathrm{N}}(T)\|$ & $2.651\times 10^{-1}$ & $15.28\times 10^{-2}$ & $11.16\times 10^{-2}$ \\
\hline\\[-3mm]
 $\|h_{\mathrm{N}}\|_{L^2(\omega)}$ & $2.128$ & $17.2777$ & $96.8253$ \\
\hline\\[-3mm]
\end{tabular}
\end{table}

\subsection{Dynamic case}

Next, we plot the uncontrolled and the controlled solutions.
\begin{figure}[H]
\centering
\begin{minipage}{0.45\textwidth}
\includegraphics[scale=0.5]{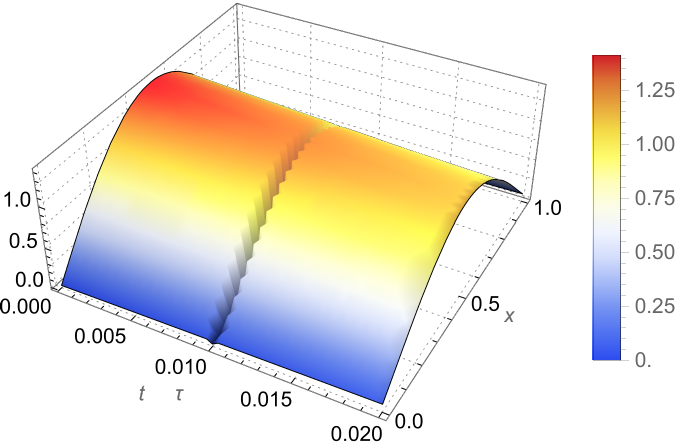}
\caption{The uncontrolled solution in Dynamic case.}
\end{minipage}\hfill
\begin{minipage}{0.45\textwidth}
\centering
\includegraphics[scale=0.5]{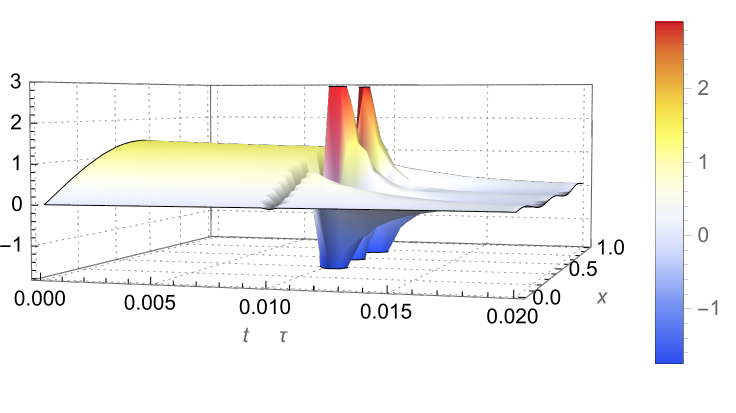}
\caption{The controlled solution in Dynamic case.}
\label{fig1dyn}
\end{minipage}
\end{figure}
The algorithm stops at the iteration number $k_*=11$. 

\begin{figure}[H]
\centering
\begin{minipage}{0.45\textwidth}
\includegraphics[scale=0.4]{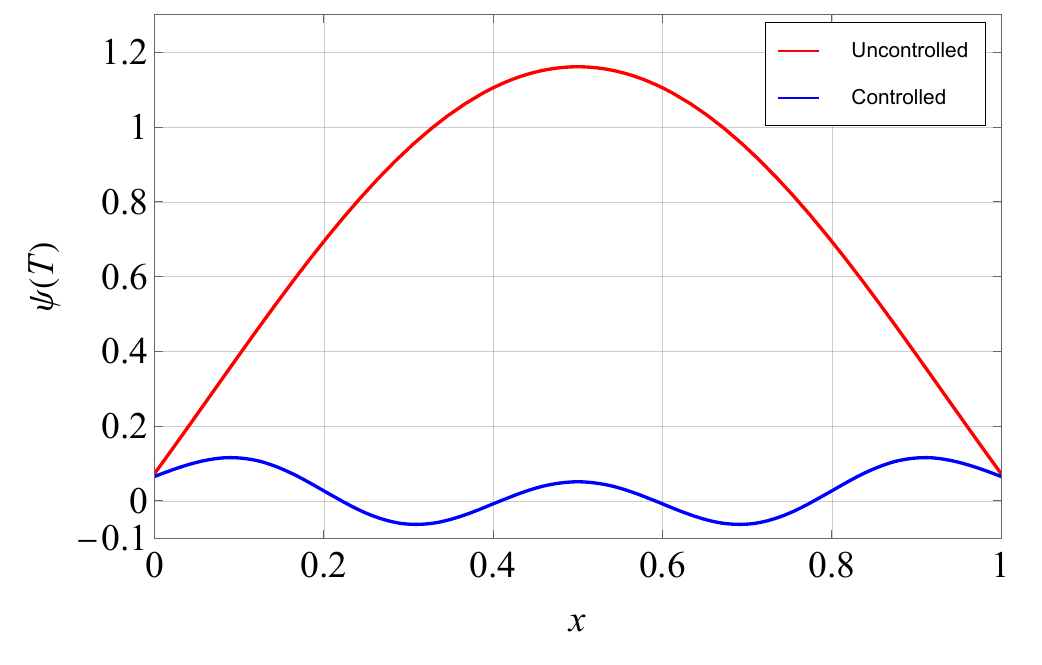}
\caption{The final state for uncontrolled and controlled solutions in Dynamic case.}
\end{minipage}\hfill
\begin{minipage}{0.45\textwidth}
\centering
\includegraphics[scale=0.4]{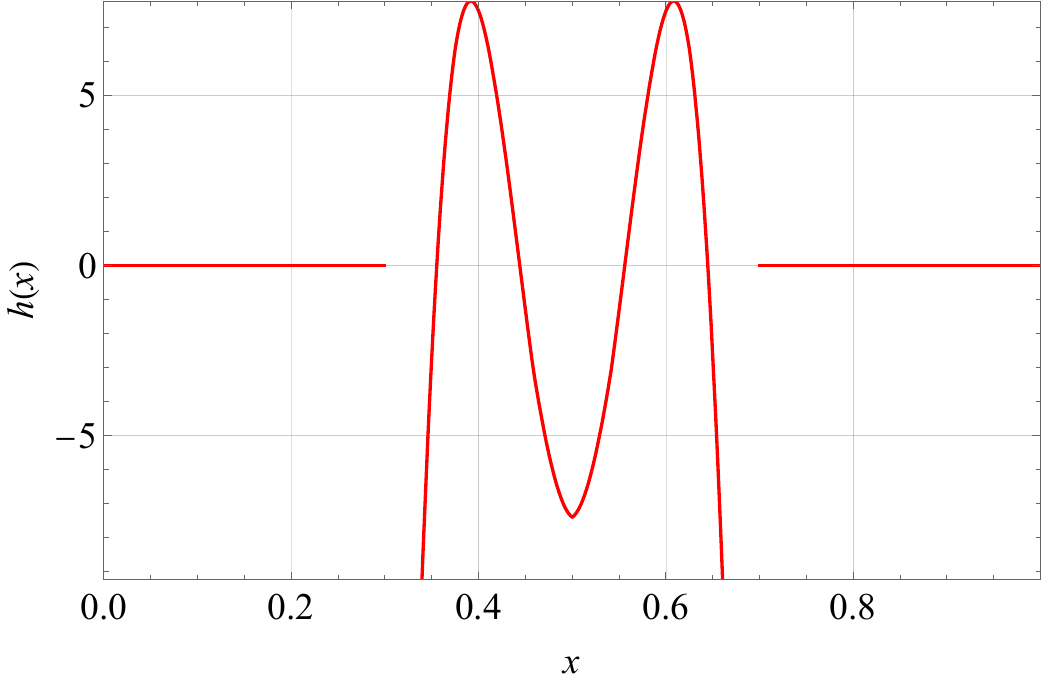}
\caption{The computed control $h$ in Dynamic case.}
\end{minipage}
\end{figure}
\bigskip

\begin{table}[ht]
\caption{Numerical results for $T=0.02$, $\tau=0.01$ and $tol=10^{-3}$ for Dynamic condition.}
\label{t3}
\centering\setlength\arraycolsep{13pt}
\begin{tabular}{cccc}
\hline\\[-3mm]
 $\varepsilon$  & $10^{-1}$ & $10^{-2}$ & $10^{-3}$ \\
\hline\\[-3mm]
 $N_{\text{iter}}$ & $4$ & $11$ & $65$ \\
\hline\\[-3mm]
 $\|\psi_{\mathrm{Dyn}} (T)\|$ & $1.598\times 10^{-1}$ & $1.135\times 10^{-1}$ & $9.02\times 10^{-2}$ \\
\hline\\[-3mm]
 $\|h_{\mathrm{Dyn}}\|_{L^2(\omega)}$ & $1.7491$ & $9.234$ & $63.6215$ \\
\hline\\[-3mm]
\end{tabular}
\end{table}

By analyzing the previous experiments, some comments and remarks are in order:
\begin{itemize}
    \item From Figures \ref{fig1dc}, \ref{fig1nc} and \ref{fig1dyn}, we notice the impact of the impulse controls at time $\tau=0.01$ on the state.
    \item Tables \ref{t1}, \ref{t2} and \ref{t3} show that when we fix the value of the penalization parameter $\varepsilon$, the Dirichlet case requires fewer iterations. The Dynamic case comes afterward with more needed iterations than the Dirichlet case. In contrast, the Neumann case requires more iteration than both previous cases.
    \item The tables also show that the norms $\|\Psi(T)\|$ decrease and the norms of the impulse controls $\|h\|_{L^2(\omega)}$ increase as $\varepsilon$ tends to zero. Moreover, for a fixed $\varepsilon$, we have
    $$\|\Psi_{\mathrm{D}}(T)\|<\|\Psi_{\mathrm{Dyn}}(T)\|<\|\Psi_{\mathrm{N}}(T)\|,$$
    and
    $$\|h_{\mathrm{D}}\|_{L^2(\omega)}<\|h_{\mathrm{Dyn}}\|_{L^2(\omega)}<\|h_{\mathrm{N}}\|_{L^2(\omega)}.$$
    These are relevant numerical observations that deserve further theoretical investigation to better understand why the above comparison holds.
\end{itemize}

The numerical simulations show that the HUM algorithm yields accurate results for the numerical approximation of impulse controls at one single instant $\tau$ for the heat equation with static boundary conditions (Dirichlet and Neumann) and also with dynamic boundary conditions. The developed algorithm deserves more investigation in the context of discrete systems and their convergence analysis in terms of discrete impulse controls. This will be investigated in future research.


\end{document}